\documentclass[a4paper,10pt]{amsart}

\usepackage{amsthm,amsmath,amsfonts,graphicx,amssymb,dsfont}
\usepackage{fullpage}
\newtheorem{defi}{Definition}[section]
\newtheorem{rmq}[defi]{Remark}
\newtheorem{theo}[defi]{Theorem}
\newtheorem*{assum}{Assumption}
\newtheorem{lem}[defi]{Lemme}
\newtheorem{prop}[defi]{Proposition}
\title{Finite mixture regression: a sparse variable selection by model selection for clustering}
\author{Emilie Devijver}
\address{Inria Select, Université Paris Sud, Bât. 425, 91405 Orsay Cedex, France}
\email{emilie.devijver@math.u-psud.fr}
\date{\today}

\begin{document}

\begin{abstract}
We consider a finite mixture of Gaussian regression model for high-dimensional data, where the number of covariates may be much larger than the sample size.
We propose to estimate the unknown conditional mixture density by a maximum likelihood estimator, restricted on relevant variables selected by an $\ell_1$-penalized maximum likelihood estimator.
We get an oracle inequality satisfied by this estimator with a Jensen-Kullback-Leibler type loss.
Our oracle inequality is deduced from a general model selection theorem for
maximum likelihood estimators with a random model collection.
We can derive the penalty shape of the criterion, which depends on the complexity of the random
model collection.

\end{abstract}
\maketitle

% \begin{keyword}[class=MSC]
% \kwd[Primary ]{}
% \kwd{}
% \kwd[; secondary ]{}
% \end{keyword}

% history:
% \received{\smonth{1} \syear{0000}}

%\tableofcontents

%\end{frontmatter}

\section{Introduction}

With the increasing of high-dimensional data, even if the number of observations is not large, new methods in statistics have been needed to deal with the identifiability underlying problem.
A classical assumption is the sparsity: if the number of parameters to estimate is larger than the sample size, we will assume that a few of parameters are nonzero.
The Lasso estimator, introduced by Tibshirani in \cite{Tibshirani}, is a classical tool in this context.
Working well in practice, many efforts have been made recently on this estimator to have some theoretical results.
Define the model and the estimator before enunce some theoretical results aready get.
We consider a linear model, $Y=X \beta + \epsilon$, with 
random variables $(X,Y) \in \mathbb{R}^p \times \mathbb{R}^q$, a regression matrix $\beta$ unknown to estimate, and a white noise $\epsilon \sim N(0,\Sigma)$.
The dimensions $p$ and $q$ could be large.
We observe the sample $((X_i,Y_i))_{i \in \{1,\ldots,n\}}$.
The Lasso estimator is defined by 
$$(\hat{\beta}^{\text{Lasso}}_{\lambda}, \hat{\Sigma}^{\text{Lasso}}_{\lambda}) = \underset{(\beta,\Sigma)}{ \operatorname{argmin}} \left\{ -\frac{1}{2n}||Y-\beta X||_2^2 + \lambda ||\beta||_1 \right\}$$
with $\lambda >0$ to specify.

Under a variety of different assumptions on the design matrix, we could have oracle inequalities for the Lasso estimator. 
For example, we can state the restricted eigenvalue condition, introduced by Bickel, Ritov and Tsybakov in \cite{Bickel}.
\begin{assum}{$RE(s,c_0)$}
 For some integer $s$ such that $1 \leq s \leq M$ and a positive number $c_0$, the following condition holds:
 $$\kappa(s,c_0)= \min_{\stackrel{J_0 \subseteq\{1,\ldots,M\}}{|J_0| \leq s}} \min_{\stackrel{\delta\neq0}{ |\delta_{J_0^c}|_1 \leq c_0 |\delta_{J_0}|_1}} \frac{|X \delta|_2}{\sqrt{n}|\delta_{J_0}|_2} >0$$
\end{assum}
With this assumption, they get an oracle inequality, which show that the distance between the prediction losses of the Lasso estimators is of the same order as the distance between it and its oracle approximation.
For an overview of existing results, cite for example  \cite{BuhlVdG} which present various conditions and various consequences.

Another type of results is about the variable selection. 
Whereas focus on the estimation, the Lasso could be used to select variables, and, for this goal, many results without hard assumptions are proved.
The first result in this way is from Meinshausen and Buhlmann, in \cite{Meins}, who show that, for neighbordhood selection in Gaussian graphical models, under a neighborhood stability condition, the Lasso is consistent,
even if the number of variables is of larger order than the sample size.
Different assumptions, as the irrepresentable Condition, described in \cite{ZhaoYu}, are in the same idea: the true variables are selected consistently.

Another approach consists to refit the estimation, after the variable selection, with an estimator with better properties.
This is the way consider in this article: we study the maximum likelihood estimator on the estimated active set.
We could cite Massart and Meynet, \cite{MassMey}, or Belloni and Chernozhukov, \cite{BellChern}, or also Tingni Sun and Cun-Hui Zhang, \cite{SunZhang} to use this idea.
Nevertheless, we will study this estimator in a finite mixture regression model, in a final goal of clustering, which is, at our knowledge, not already studied.

The goal of clustering methods is to discover structures among individuals described by several variables.
Specifically, in regression case, given $n$ observations $(x,y)=((x_1,y_1),\ldots,(x_n,y_n))$ which are realizations of random variables $(X,Y)$ with $X \in \mathbb{R}^p$ and $Y \in \mathbb{R}^q$, one aims at grouping the data into a few clusters such that the conditional observations $Y|X$ in the same cluster are more similar to each other
than those from the other clusters.
Different methods could be envisaged, more geometric or more statistical.
We are dealing with model-based clustering, in order to have a rigorous statistical framework to assess the number of clusters and the role of each variable.
Datasets are more and more in high-dimension, and all the information should not be interesting for the clustering.
To solve this problem, we propose a procedure which provide a data clustering from variable selection. This procedure is based on a modeling that recasts variable selection and clustering problems into a model selection problem in a regression framework.
A global selection criterion choosing simultaneously the best number of clusters and the set of relevant variables is required.
We use a penalized criterion to select a model from a non-asymptotic point of view.
Penalizing the empirical contrast is an idea emerging from the seventies. Akaike, in \cite{akaike}, proposed the Akaike's Information
Criterion (AIC) in $1973$, and Schwarz in $1978$ in \cite{schwarz} suggested the Bayesian Information Criterion (BIC). Those criteria are based on asymptotic heuristics.
To deal with non-asymptotic observations, Birgé and Massart in \cite{BirgeMassart} and Barron et al. in \cite{Barron}, define a penalized data-driven criterion,
which leads to oracle inequalities for model selection.
Cohen and Le Pennec, in \cite{Lepennec}, generalize this result in the case of regression data.
The aim of our approach is to define penalized data-driven criterion which leads to an oracle inequality for our procedure.
In our context of regression, Cohen and Le Pennec, in \cite{Lepennec}, proposed a general model selection theorem for maximum likelihood estimation, adapted from Massart's theorem in \cite{MassartStFlour}.
Nevertheless, we can not apply it directly, because it is stated for a deterministic model collection, whereas our data-driven model collection is random, constructed by the Lasso.
As Meynet done in \cite{MeynetMaugis} to generalize Massart's theorem, we extend the theorem to cope with the randomness of our model collection.
By applying this general theorem to the finite mixture regression random model collection constructed by our procedure, we derive a convenient theoretical penalty as well as an associated non-asymptotic penalized criteria and an oracle inequality fulfilled by our Lasso-MLE estimator.
The advantage of this procedure is that it does not need any restrictive assumption.

Let give the main result of this paper.
Let $(x_i,y_i)_{i=1,\ldots,n}$ the observations, with unknown conditional density $s_0$.
Let $(S_m)_{m\in \mathcal{M}}$ the model collection constructed by our procedure. 
We construct a collection of finite regression mixture of Gaussians with various numbers of clusters and different sets of relevant variables.
Then, we estimate the conditional density by the maximum likelihood estimator in each model. This leads to a collection of estimators for the density. A final estimator has to be selected among this collection, which is equivalent to select a model among the model collection.
Under some weak assumptions, we obtain a minimizer of $pen(m)$ such that the estimator $\hat{s}_{\hat{m}}$, $\hat{s}$ being the maximum likelihood estimator,
and $\hat{m}$ the model which minimizes the penalized log-likelihood,
satisfies
\begin{align*}
&E \left[JKL_{\rho,\lambda}^{\otimes_n} (s_0,\hat{s}_{(\hat{k},\hat{J})}) \right] \\
\leq &C \left[ E \left( \inf_{(k,J) \in \hat{\mathcal{M}}} \left(  \inf_{t \in S_{(k,J)}} KL_{\lambda}^{\otimes_n} (s_0,t) + \frac{\text{pen}(k,J)}{n} \right)\right) + \frac{4}{n} \right]. 
\end{align*}

We will define $JKL$ and $KL$ later.
The idea of this theorem is that the model choose by our procedure is as good as the best we can do among our collection, even if we have known the true density.

Before concluding the introduction, let give some notations which need to be fixed.
In this general setting, we assume that the observations $(x_i,y_i)_{i=1,\ldots,n}$ are a sample of random variables $(X,Y)$ where $X \in \mathcal{X}$ and $Y \in \mathcal{Y}$.
Let  $S_m$ a set of candidate conditional densities, in which we estimate $\hat{s}_m$ with the maximum likelihood estimator
$$\hat{s}_m = \underset{s_m \in S_m}{\operatorname{argmin}} \left( - \sum_{i=1}^n \log s_m(y_i|x_i) \right).$$
To avoid existence issue, we work with almost minimizer of this quantity and define an $\eta$-log-likelihood minimizer as any $\hat{s}_m$ that satisfies
$$\sum_{i=1}^n - \log(\hat{s}_m (y_i|x_i)) \leq \inf_{s_m \in S_m} \left( \sum_{i=1}^n - \log(s_m(y_i|x_i)) \right) + \eta.$$
The best model in this collection is the one with the smallest risk. However, because we do not have access to the true density $s_0$, we can not select the best model, which will be called the oracle. Thereby,
there is a trade-off between a bias term measuring the closeness of $s_0$ to the set $S_m$ and a variance term depending on the complexity of the set $S_m$ and on the sample size. 
A good set $S_m$ will be thus one for which this trade-off leads to a small risk bound. 
We are working with a maximum likelihood approach, the most natural quality measure is thus the Kullback-Leibler divergence denoted by $KL$. As we consider law with densities with respect to the Lebesgue measure $d\lambda$, we use the following notation 
\begin{align*}
KL_\lambda (s,t) &= KL(sd\lambda,td\lambda)\\
&= \left\{ \begin{array}{ll}
    &\int \log \left(\frac{s}{t}\right) s d\lambda \text{ if } sd\lambda <<t d\lambda ;\\
    &+\infty \text{ otherwise. }
   \end{array} \right.
\end{align*}

Remark that, contrary to the quadratic loss, this divergence is an intrinsic quality measure between probability laws: it does not depend on the reference measure $d\lambda$. 
However, the densities depend on this reference measure, and this is stressed by the index $\lambda$.
As we deal with conditional densities and not classical densities, the previous divergence should be adapted.

We define the tensorized Kullback-Leibler divergence by
$$ KL^{\otimes_n} _{\lambda} (s,t) = E\left[ \frac{1}{n} \sum_{i=1}^n KL_{\lambda} (s(.|x_i), t(.|x_i)) \right].$$

This divergence used in \cite{Lepennec} appears as the natural one in this regression setting.

Namely, we use the Jensen-Kullback-Leibler divergence $JKL_{\rho}$ with $\rho \in ]0,1[$ defined by
$$JKL_{\rho} (s d\lambda, t d\lambda) = JKL_{\rho,\lambda} (s,t) = \frac{1}{\rho} KL_{\lambda}(s,(1-\rho)s + \rho t) ;$$
and the tensorized one
$$JKL^{\otimes_n}_{\rho,\lambda} (s,t) = E \left[ \frac{1}{n} \sum_{i=1}^ n JKL_{\rho,\lambda}^ {\otimes_n} (s(.|x_i),t(.|x_i)) \right].$$
This divergence is studied in \cite{Lepennec}.
We prefer this divergence rather than the Kullback-Leibler one because we get a boundness assumption on the controlled functions
that is not satisfied by the log-likelihood differences differences $-\log \frac{s_m}{s_0}$. When considering the Jensen-Kullback-Leibler divergence, those ratios are replaced by ratios $-\frac{1}{\rho} \log \left( \frac{(1-\rho) s_0 + \rho s_m}{s_0} \right)$
that are close to the log-likelihood differences when the $s_m$ are close to $s_0$ and always upper bounded by $- \frac{\log(1-\rho)}{\rho}$.

Indeed, it is needed to use deviation inequalities
for sums of random variables and their suprema, which is the key of the proof of oracle type inequality.

The aim of the model selection is to construct a data-driven criterion to select a model of proper dimension of a given list. A general theory of this topic is proposed in the works of Birgé and Massart \cite{BirgeMassart2}.
Besides, Massart, in \cite{MassartStFlour}, proposed a general theorem which gives the form of the penalty and associated oracle inequality in term of the Kullback-Leibler and Hellinger loss.
In our case of regression, Cohen and Le Pennec, in \cite{Lepennec}, proposed a general theorem which gives the form of the penalty and associated oracle inequality in term of the Kullback-Leibler and Jensen-Kullback-Leibler loss.
These theorems are based on the centred process control with the bracketing entropy, allowing to evaluate the size of the models.
We compare the risk of the penalized maximum likelihood estimator $\hat{s}_{\hat{m}}$ with the benchmark $\inf_{m \in \hat{\mathcal{M}}} E(KL_{\lambda}^{\otimes_n} (s,\hat{s}_m))$.
Our setting is more general, because we work with a random family denoted by $\hat{\mathcal{M}}$. We have to control the centred process thanks to Bernstein's inequality.

The rest of the article is organized as follows.
In the section \ref{sectionProc}, we recall the multivariate Gaussian mixture regression model, and we describe the main steps of the procedure we propose. We also illustrate the requirement of refitting by some simulations.
We present our oracle inequality in the section \ref{inegOracl}. 
Finally, in section \ref{sectionProof}, we give some tools to understand the proof of the oracle inequality, with a global theorem of model selection with a random collection in section \ref{sectionThmGeneral} and sketch of proofs after.
All the details are given in Appendix.

\section{The Lasso-MLE procedure}
\label{sectionProc}
In order to cluster high-dimensional regression data, we will work with the multivariate Gaussian mixture regression model.
This model is developed in \cite{VandeGeer} in the scalar response case.
We generalize it in section \ref{model}.
Moreover, we want to construct a model collection. We propose, in section \ref{proc}, a procedure called Lasso-MLE which construct a model collection, with various sparsity, 
of Gaussian mixture regression models.
The different sparsities solve the high-dimensional problem.
We conclude this section with a simulation, which illustrate the advantage of refitting.
\subsection{Gaussian mixture regression model}
\label{model}
We observe $n$ independent couples $(x_{i}, y_{i})_{1 \leq i \leq n}$ of random variables $(X,Y)$, with $Y  \in \mathbb{R}^q$ and $X \in \mathbb{R}^p$ coming from a probability distribution with unknown conditional density denoted by $s_0$.
To solve a clustering problem, we use a finite mixture model in regression. In particular, we will approximate the density of $Y|X$ with a multivariate Gaussian mixture regression model. If the observation $i$ belongs to the cluster $r$, we assume that there exists $\beta_r \in \mathbb{R}^{p \times q}$ such that $y_i= \beta_r x_i +\epsilon$ where $\epsilon \sim N(0,\Sigma_r)$.

Thus, the random response variable $ Y \in \mathbb{R}^q$ depends on a set of explanatory variables, written $X \in \mathbb{R}^p$, through a regression-type model.
 Give more precisions on the assumptions.
\begin{itemize}
 \item The variables $Y_i|X_i$ are independent, for all $i=1,\ldots, n$ ;
 \item the variables $Y_i|X_i=x_i \sim s_{\xi}(y|x_i)dy$, with
  \begin{align}
  \label{modele}
  &s_{\xi}(y|x)=\sum_{r=1}^{k} \frac{\pi_{r}}{(2 \pi)^{\frac{q}{2}}  \text{det}(\Sigma_r)^{1/2}} \exp \left( -\frac{(y-\beta_{r} x)^t \Sigma_{r}^{-1}(y-\beta_{r} x)}{2} \right)\\
  &\xi=(\pi_{1},\ldots, \pi_{k},\beta_{1},\ldots,\beta_{k},\Sigma_{1},\ldots,\Sigma_{k}) \in \left( \Pi_{k} \times (\mathbb{R}^{q\times p})^k\times (\mathbb{S}^q_{++})^k \right) \nonumber\\
  & \Pi_{k} = \left\{ (\pi_{1}, \ldots, \pi_{k}) ; \pi_{r} >0 \text{ for } r \in \{1,\ldots, k\} \text{ and } \sum_{r=1}^{k} \pi_{r} = 1 \right\} \nonumber \\
  & \mathbb{S}_{++}^q \text{ is the set of symmetric positive definite matrices on } \mathbb{R}^q . \nonumber
 \end{align}
\end{itemize}

We want to estimate the conditional density function $s_\xi$ from the observations. For all $r \in \{1,\ldots, k\}, \beta_{r}$ is the matrix of regression coefficients, and $\Sigma_{r}$ is the covariance matrix in the mixture component $r$. The $\pi_{r}$s are the mixture proportions.
In fact, for all $r \in \{1,\ldots,k\}$, for all $z \in \{1,\ldots,q\}$, $\beta_{r,z}^t x= \sum_{j=1}^{p} \beta_{r,j,z} x_{j}$ is the $z$th component of the mean of the mixture component $r$ for the conditional density $s_{\xi}(.|x)$.

A variable is said to be irrelevant if, for each $r \in \{1,\ldots,k\}$, $\beta_r=0$. A variable is relevant if it is not irrelevant.
A model is said to be sparse if there is a few of relevant variables.

 We denote by $x^{[J]}$ the restriction of $x$ on $J$, and 
$\mathcal{S}_{(k,J)}$ the model with $k$ components and with $J$ for relevant variables set:

\begin{equation}
\label{modele h}
\mathcal{S}_{(k,J)} = \left\{  y \in \mathbb{R}^q | x \in \mathbb{R}^p \mapsto s^{(k,J)}_{\xi}(y|x) \right\}
 \end{equation}
where 
$$s^{(k,J)}_{\xi}(y|x)=\sum_{r=1}^{k} \frac{\pi_{r}}{(2 \pi)^{\frac{q}{2}}  \text{det}(\Sigma_r)^{1/2}} \exp \left( -\frac{(y-(\beta_{r} x)_{|J})^t \Sigma_{r}^{-1}(y-(\beta_{r} x)_{|J})}{2} \right)$$

This is the main model used in this paper.
Nevertheless, to deal with high-dimensional data, we use the Lasso estimator to construct the set of  relevant variables, and the choice of the regularization parameter is known to be a difficult problem.
We propose to construct a model collection to solve this problem.
 \subsection{The Lasso-MLE procedure}
 \label{proc}
The procedure we propose which is particularly interesting in high-dimension could be decomposed into three main steps.

The first step consists of constructing a collection of models $\{\mathcal{S}_{(k,J)}\}_{(k,J) \in \mathcal{M}}$ in which the model $\mathcal{S}_{(k,J)}$ is defined by equation \eqref{modele h},
and the model collection is indexed by $ \mathcal{M}= K \times \mathcal{J}$. Denote $K\subset \mathbb{N}^*$ the possible number of components, and
denote $\mathcal{J}$ a collection of subsets of $\{1,\ldots,p\} \times \{1,\ldots,q\}$.

To detect the relevant variables, and construct the set $J$ in each model, we penalize the empirical contrast by an $\ell_1$-penalty on the mean parameters proportional to $||P_r\beta_r||_1 = \sum_{j=1}^p \sum_{z=1}^q |(P_r \beta_r)_{j,z}|$, where $P_r^t P_r = \Sigma_r^{-1}$. 
This leads to penalize simultaneously the $\ell_1$-norm of the mean coefficients and small variances.
Computing those estimators lead to the relevant variables set.
For a fixed number of mixture components $k \in K$, denote by $G_k$ a candidate of regularization parameters. Fix a parameter $\lambda \in G_k$, we could then use an EM algorithm to compute the set of relevant variables.
Then, varying $k \in K$ and $\lambda \in G_k$, we construct the relevant variables set $J_{k,\lambda}$. We denote by $\mathcal{J}$ the random collection of all these sets, $\mathcal{J}= \bigcup_{k \in K} \bigcup_{\lambda \in G_k} J_{(k,\lambda)}$ .

The second step consists of approximating the MLE 
$$\hat{s}_{(k,J)}= \underset{t \in \mathcal{S}_{(k,J)}}{ \operatorname{argmin}} \left\{ \frac{1}{n} \sum_{i=1}^n \log (t(y_i|x_i)) \right\} $$ 
using an EM algorithm for each model $(k,J)\in \mathcal{M}$. 

The third step is devoted to model selection. 
We get a model collection, and we need to choose the best one. Because we do not have access to $s_0$, we can not take the one which minimizes the risk.
The theorem \ref{thm general} solve this problem: we get a penalty achieving to an oracle inequality. 
Then, even if we do not have access to $s_0$, we know that we can do almost like the oracle.
%In practice, we do not use the theoretical penalty. We use the slope heuristic to calibrate the criterion to minimize.
\subsection{Why refit the Lasso estimator?}

In order to illustrate our procedure, we compute multivariate data, the restricted eigenvalue condition being not satisfied, and run our procedure.
We consider an extension of the model studied in Giraud et al. article \cite{Giraud} in the section $6.3$.
Indeed, this model is a linear regression with a scalar response which does not satisfy the restricted eigenvalues condition.
Then, we define different classes, to get a finite mixture regression model, which does not satisfied the restricted eigenvalues condition, and extend the dimension for multivariate response.
We could compare the result of our procedure with the Lasso, to illustrate the oracle inequality we have get.
Let precise the model.

Let $x^{(1)}, x^{(2)}, x^{(3)}$ be three vectors of $\mathbb{R}^n$ defined by
$$\begin{array}{ll}
 x^{(1)} &= ( 1,-1,0,\ldots, 0)^t /\sqrt{2} \\
 x^{(2)} &= ( -1,1.001,0,\ldots, 0)^t /\sqrt{1+0.001^2} \\
 x^{(3)} &= ( 1/\sqrt{2},1/\sqrt{2},1/n,\ldots, 1/n)^t /\sqrt{1+(n-2)/n^2} \\ 
\end{array}
$$
and for $4\leq j \leq n$, let $x^{(j)}$ be the $j^{th}$ vector of the canonical basis of $\mathbb{R}^n$.
We take a sample of size $n=20$, and vector of size $p=m=10$. We consider two classes, each of them define by $\beta_{j,z,1}=10$ and $\beta_{j,z,2}=-10$ for $j \in \{1,\ldots,2\}$, $z\in \{1,\ldots,10\}$.
Moreover, we define the variance of the noise by a diagonal matrix with $0.01$ for diagonal coefficient in each class.

\begin{figure}
 \centering
 \includegraphics[scale=0.4]{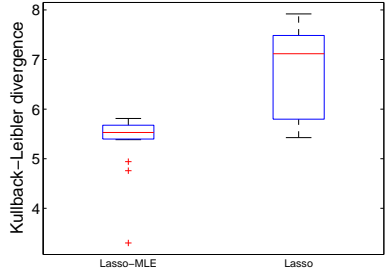}
 \caption{Boxplot of the Kullback-Leibler divergence between the true model and the one constructed by each procedure, the Lasso-MLE procedure and the Lasso procedure.}
\label{KL}
 \end{figure}

We run our procedure on this model, and compare it with the Lasso procedure, without refitting. We compute the model selected by the slope heuristic over the model collection constructed by the Lasso estimator.
In figure \ref{KL} are the boxplots of each procedure, running $20$ times. The Kullback-Leibler divergence is computed over a sample of size $5000$.

 We could see that a refitting after variable selection by the Lasso leads to a better estimation, according to the Kullback-Leibler loss.
\section{An oracle inequality for the Lasso-MLE estimator}
\label{inegOracl}
Let denote the model collection constructed by the Lasso-MLE procedure by $\mathcal{S}=(\mathcal{S}_{(k,J)})_{(k,J) \in \mathcal{M}^L}$.
The model
$S_{(k,J)}$ is defined in \eqref{modele h},
whereas we have denoted $\mathcal{M}^L= K\times \mathcal{J}^L$, with $\mathcal{J}^L$ a random subcollection of $\mathcal{P}(\{1,\ldots,p\} \times \{1,\ldots,q\})$, constructed by the Lasso.

We will work with restricted parameters. 
Assume $\Sigma_r$ diagonal, with $\Sigma_r=\text{diag} (\Sigma_{1,r}^2, \ldots, \Sigma_{q,r}^2)$, for all $r \in \{1,\ldots,k\}$. We define

\begin{align}
\label{ens borne}
\mathcal{S}_{(k,J)}^\mathcal{B} = &\left\{  s_\xi^{(k,J)} \in \mathcal{S}_{(k,J)}, (\beta_r)_{|J} \in [-A_\beta,A_\beta]^{J}, \right. \nonumber \\
&\left. a_\Sigma^2 \leq \Sigma_{z,r} \leq A_\Sigma^2  \text{ for all } z \in [1,q] \text{ for all } r \in [1,k] \right\}. 
\end{align}
Moreover, we assume that the covariates $X$ belong to an hypercube. Without restriction, we could assume that $X \in [0,1]^p$.
\begin{rmq}
 We have to denote that in this paper, the active variables set is designed by the Lasso.
 Nevertheless, any tool is used to construct this set, we could obtain analog results. We could work with any random subcollection of $\mathcal{P}(\{1,\ldots,p \} \times \{1,\ldots,q\})$, the control ed size being required in high-dimensional case.
\end{rmq}

\begin{theo}
\label{thm}
Let $(x_i,y_i)_{i=1,\ldots,n}$ the observations, with unknown conditional density $s_0$.
Let $\mathcal{S}_{(k,J)}$ as defined in \eqref{modele h}.
We denote by $\mathcal{M}^L$ a random subcollection of $\mathcal{M}$.
For $(k,J) \in \mathcal{M}^L$, denote $\mathcal{S}_{(k,J)}^\mathcal{B}$ the model defined in \eqref{ens borne}.

Consider the maximum likelihood estimator
$$\hat{s}_{(k,J)} = \underset{s_\xi \in \mathcal{S}_{(k,J)}^\mathcal{B}}{\operatorname{argmin}}\left\{ -\frac{1}{n} \sum_{i=1}^n \log (s_\xi(y_i|x_i)) \right\}.$$
Denote by $D_{(k,J)}$ the dimension of the model $\mathcal{S}_{(k,J)}^\mathcal{B}$, $D_{(k,J)} = k(|J|+q+1)-1$. Let $\bar{s} \in \mathcal{S}_{(k,J)}^\mathcal{B}$ such that
$$KL^{\otimes_n}_{\lambda} (s_0,\bar{s}) \leq \inf_{t \in \mathcal{S}_{(k,J)}^\mathcal{B}} KL^{\otimes_n}_{\lambda} (s_0,t) +\frac{\delta_{KL}}{n} ;$$
and let $\tau>0$ such that $\bar{s}\geq e^{-\tau} s_0$.
Let $\text{pen}: \mathcal{M} \rightarrow \mathbb{R}_+$, and suppose that there exists an absolute constant $\kappa>0$ such that, for all $(k,J) \in \mathcal{M}$,
\begin{align*}
\text{pen}(k,J) \geq \kappa D_{(k,J)} & \left[ B^2(A_\beta,A_\Sigma,a_\sigma,q) -\log\left(\frac{D_{(k,J)}}{n} B^2(A_\beta,A_\Sigma,a_\sigma,q) \wedge 1\right) \right.\\
& \left. + (1\vee \tau) \log \left( \frac{4epq}{(D_{(k,J)}-q^2) \wedge pq }\right) \right].
\end{align*}

Then, the estimator $\hat{s}_{(\hat{k},\hat{J})}$, with
$$(\hat{k},\hat{J}) =  \underset{(k,J) \in \mathcal{M}^L}{\operatorname{argmin}} \left \{ -\frac{1}{n} \sum_{i=1}^n \log (\hat{s}_{(k,J)} (y_i |x_i)) + \text{pen}(k,J) \right\}$$

satisfies
\begin{align*}
E \left[JKL_{\rho,\lambda}^{\otimes_n} (s_0,\hat{s}_{\hat{m}}) \right] \leq & C_1 E \left( \inf_{(k,J) \in \mathcal{M}^L} \left(  \inf_{t \in S_{(k,J)}} KL_{\lambda}^{\otimes_n} (s_0,t) + \frac{\text{pen}(k,J)}{n} \right)\right) \\
&+ C_2 \frac{\Sigma^2}{n} 
\end{align*}

for some absolute positive constants $C_1$ and $C_2$.

\end{theo}

This result could be compare with the oracle inequality get in \cite{VandeGeer}.
Indeed, under restricted eigenvalues condition (this assumption is explained in details in Bühlman and Van de Geer's book \cite{VdGBook}) and fix design, they get an oracle inequality for the Lasso estimator in finite mixture regression model, with scalar response and high-dimension regressors.
We get a similar result for the Lasso-MLE estimator.
The good point is that we get the same type of inequality as comparable estimators.
Moreover, our procedure work in a more general framework, without any assumptions about the design.

\section{Tools for proof}
In this section, we present the tools needed to understand the proof. First, we present a general theorem for model selection in regression among a random collection.
Then, in subsection \ref{proof general theorem}, we present the proof of this theorem, and in the next subsection we explain how use the main theorem to get the oracle inequality. All details are available in Appendix.
\label{sectionProof}
\subsection{General theory of model selection with the maximum likelihood estimator.}
\label{sectionThmGeneral}
To get an oracle inequality for our clustering procedure, we have to use a general model selection theorem. Because the model collection
constructed by our procedure is random, because of the Lasso estimator which select variables randomly, we have to generalize theorems already existing.
Begin by some general theory of model selection.

Before enunciate the general theorem, begin by talking about the assumptions.
First, we impose a structural assumption.
It is a bracketing entropy condition on the model $S_m$ with respect to the Hellinger divergence $d^{2 \otimes_n}_{H}(s,t) = E \left[ \frac{1}{n} \sum_{i=1}^n d^2_{H} (s(.|x_i),t(.|x_i)) \right]$.
A bracket $[t^-,t^+]$ is a pair of functions such that for all $(x,y) \in \mathcal{X} \times \mathcal{Y}, t^-(y,x) \leq s(y|x) \leq t^+ (y,x)$.
The bracketing entropy $H_{[.]} (\delta,S,d_H^{\otimes_n})$ of a set $S$ is defined as the logarithm of the minimum number of brackets $[t^-,t^+]$ of width $d_H^{\otimes_n}(t^-,t^+)$ smaller than $\delta$ such that every functions of $S$ belong to one of these brackets.

\begin{assum}[$\text{H}_m$]
There is a non-decreasing function $\phi_m$ such that $\delta \mapsto \frac{1}{\delta} \phi_m (\delta)$ is non-increasing on $(0,+\infty)$ and for every $\sigma \in \mathbb{R}^+$ and every $s_m \in S_m$,
$$\int _0^\sigma \sqrt{ H_{[.]} (\delta,S_m(s_m,\sigma),d_H^{\otimes_n})} d \delta \leq \phi_{m}(\sigma)$$
where $S_m(s_m,\sigma)=\{t \in S_m, d_H^{\otimes_n}(t,s_m) \leq \sigma \}$.
The model complexity $\mathcal{D}_m$ is then defined as $n \sigma^2_m$ with $\sigma^2_m$ the unique root of 
\begin{align}
\label{model complexity}
\frac{1}{\sigma} \phi_{m}(\sigma) = \sqrt{n}\sigma.
\end{align}
\end{assum}
Denote that the model complexity depends on the bracketing entropies not of the global models $S_m$ but of the ones of smaller localized sets.
This is a weaker assumption.

For technical reason, a separability assumption is also required.
\begin{assum}[$\text{Sep}_m$]
There exists a countable subset $S^{'}_m$ of $S_m$ and a set $\mathcal{Y}_m^{'}$ with $\lambda(\mathcal{Y} \setminus \mathcal{Y}^{'}_m)=0$ such that for every $t \in S_m$, there exists a sequence $(t_k)_{k \geq 1}$ of elements of $S_m^{'}$ such that for every $x$ and every $y \in \mathcal{Y}_m^{'}$, $\log(t_k(y|x))$ goes to $\log(t(y|x))$ as $k$ goes to infinity. 
\end{assum}

We also need an 
information theory type assumption on our collection. We assume the existence of a Kraft-type inequality for the 
collection:
\begin{assum}[K]
There is a family $(x_m)_{m \in \mathcal{M}}$ of non-negative numbers such that
$$\sum_{m \in \mathcal{M}} e^{-x_m} \leq \Sigma < + \infty.$$
 \end{assum}
The difference with Cohen and Le Pennec's theorem is that we consider a random collection of models $\hat{\mathcal{M}}$, included in the whole collection $\mathcal{M}$.
In our procedure, we deal with the high-dimensional models, and we cannot look after all the models: we have to restrict ourselves to a smaller subcollection of models.

Then we could write our main global theorem.
\begin{theo}
\label{thm general}
 Assume we observe $(x_i,y_i)$ with unknown conditional density $s_0$. Let $\mathcal{S} = (S_m)_{m\in \mathcal{M}}$ be at most countable collection of conditional density sets.
 Assume assumption (K) holds, while assumptions $(H_m)$ and $(\text{Sep}_m)$ hold for every model $S_m \in \mathcal{S}$. 
 Let $\delta_{KL}>0$, and $\bar{s}_m \in S_m$ such that 
 $$KL_\lambda^{\otimes_n} (s_0,\bar{s}_m) \leq  \inf_{t \in S_m} KL_\lambda^{\otimes_n}(s_0,t) +\frac{\delta_{KL}}{n} ;$$
 and let $\tau >0$ such that 
 \begin{equation}
 \label{hyp tau}
 \bar{s}_m \geq e^{-\tau} s_0.
 \end{equation}

 Introduce $(S_m)_{m \in \hat{\mathcal{M}}}$ some random subcollection of $(S_m)_{m \in \mathcal{M}}$.
 Consider the collection $(\hat{s}_m) _{m\in  \hat{\mathcal{M}}}$ of $\eta$-log-likelihood minimizer in $S_m$, satisfying, for all $m \in \hat{\mathcal{M}}$,
 $$\sum_{i=1}^n - \log(\hat{s}_m (y_i|x_i)) \leq \inf_{s_m \in S_m} \left( \sum_{i=1}^n - \log(s_m(y_i|x_i)) \right) + \eta.$$
 
 Then, for any $\rho \in (0,1)$ and any $C_1 > 1$, there are two constants $\kappa_0$ and $C_2$ depending only on $\rho$ and $C_1$ such that, as soon as for every index $m \in \mathcal{M}$, 
 \begin{equation}
\label{penalite}
 \text{pen}(m) \geq \kappa (\mathcal{D}_m + (1 \vee \tau) x_m)
   \end{equation}

 with $\kappa > \kappa_0$, and where the model complexity $\mathcal{D}_m$ is defined in \eqref{model complexity},
 the penalized likelihood estimate $\hat{s}_{\hat{m}}$ with $\hat{m} \in \hat{\mathcal{M}}$ such that
 $$\sum_{i=1}^n - \log (\hat{s}_{\hat{m}} (y_i|x_i)) + \text{pen}(\hat{m}) \leq \inf_{m \in \hat{\mathcal{M}}} \left( \sum_{i=1}^n - \log( \hat{s}_m (y_i |x_i)) + \text{pen} (m) \right) + \eta^{'}$$
 satisfies
\begin{align}
\label{inegalité oracle}
 E(JKL_{\rho,\lambda}^ {\otimes_n} (s_0,\hat{s}_{\hat{m}} )) \leq & C_1 E \left( \inf_{m \in \hat{\mathcal{M}}} \inf_{t \in S_m} KL^ {\otimes_n }_{\lambda} (s_0,t)+2 \frac{\text{pen}(m)}{n} \right) \nonumber \\
 &+ C_2 (1 \vee \tau) \frac{\Sigma^2}{n} + \frac{\eta' + \eta}{n}.
\end{align}
\end{theo}
Obviously, one of the models minimizes the right hand side. Unfortunately, there is no way to know which one without knowing $s_0$. Hence, this oracle model can not be used to estimate $s_0$. 
We nevertheless propose a data-driven strategy to select an estimate among the collection of estimates $\{\hat{s}_m\}_{m\in \hat{\mathcal{M}}}$ according to a selection rule that performs almost as well as if we had known this oracle, according to the absolute constant $C_1$.
Using simply the log-likelihood of the estimate in each model as a criterion is not sufficient.
It is an underestimation of the true risk of the estimate and this leads to choose models that are too complex.
By adding an adapted penalty $\text{pen}(m)$, one hopes to compensate for both the variance term and the bias term between $\frac{1}{n} \sum_{i=1}^n - \log\frac{\hat{s}_{\hat{m}}(y_i|x_i)}{s_0(y_i|x_i)}$ and
$\inf_{s_m \in S_m} KL_{\lambda}^{\otimes_n } (s_0,s_m) $. For a given choice of $\text{pen}(m)$, the best model $S_{\hat{m}}$ is chosen as the one whose index is an almost minimizer of the penalized $\eta$-log-likelihood.

Talk about the assumption \eqref{hyp tau}.
If $s$ is bounded, with a compact support, this assumption is satisfied.
It is also satisfied in other cases, more particular. Then it is not a hard assumption, and it is needed to control the random family.

This theorem is available for whatever model collection constructed, whereas assumptions $(H_m), (K)$ and $(Sep_m)$ are satisfied.
In the following, we will specify the procedure we propose to cluster high-dimensional data, and look for satisfying these assumptions.
Nevertheless, this theorem is not specific of our context, and could be used whatever the problem considering.

\subsection{Proof of the general theorem}
\label{proof general theorem}
For the sake of simplicity, we shall assume that $\rho=0$.
For any model $S_m$, we have denoted that $\bar{s}_m$ a function such that
$$KL^{\otimes_n }_\lambda (s_0,\bar{s}_m) \leq \inf_{s_m \in S_m} KL_{\lambda}^{\otimes_n } (s_0,s_m) + \frac{\delta_{KL}}{n}.$$

Fix $m \in \mathcal{M}$ such that $KL^{\otimes_n }_\lambda (s_0,\bar{s}_m) < + \infty$. Introduce 
\begin{align*}
M(m)=& \left\{ m' \in \mathcal{M}, P_n^{\otimes_n} (- \log \hat{s}_{m^{'}}) + \frac{\text{pen}(m^{'})}{n} \right. \\
&\left. \hspace{1.3cm} \leq  P_n^{\otimes_n}(- \log \hat{s}_{m}) + \frac{\text{pen}(m)}{n} + \frac{\eta^{'}}{n} \right\} ;
\end{align*}
where $P_n^ {\otimes_n} (g) = \frac{1}{n} \sum_{i=1}^{n}g(Y_i|X_i)$. 
We define the functions $kl(\bar{s}_m),kl(\hat{s}_m)$ and $jkl(\hat{s}_m)$ by
\begin{align*}
kl(\bar{s}_m)&=-\log \left(\frac{\bar{s}_m}{s_0} \right) ;\hspace{1cm} kl(\hat{s}_m)=- \log \left( \frac{\hat{s}_m}{s_0} \right) ;\\
jkl(\hat{s}_m)&=-\frac{1}{\rho} \log \left( \frac{(1-\rho)s_0 +\rho \hat{s}_m}{s_0} \right). 
\end{align*}

For every $m^{'} \in \mathcal{M}(m)$, by definition, 
 \begin{align*}
P_n^{\otimes_n}(kl(\hat{s}_{m^{'}})) +\frac{\text{pen}(m^{'})}{n} &\leq  P_n^{\otimes_n}(kl(\hat{s}_{m})) +\frac{\text{pen}(m) + \eta^{'}}{n}\\
&\leq  P_n^{\otimes_n}(kl(\bar{s}_{m})) +\frac{\text{pen}(m) + \eta^{'}+\eta}{n}.  
 \end{align*}

Let $\nu_n^ {\otimes_n } (g)$ denote the recentred process $P_n^ {\otimes_n}(g) - P^{\otimes_n}(g)$.
By concavity of the logarithm, \mbox{$kl(\hat{s}_{m^{'}}) \geq jkl(\hat{s}_{m^{'}})$}, and then

\begin{align*}
&P^{\otimes_n}(jkl(\hat{s}_{m^{'}})) - \nu_{n}^{\otimes_n} (kl(\bar{s}_m)) \\
\leq  &P^{\otimes_n}(kl(\bar{s}_{m})) +\frac{\text{pen}(m) }{n}- \nu_{n}^{\otimes_n} (jkl(\hat{s}_{m^{'}})) +\frac{\eta^{'}+\eta}{n} -\frac{\text{pen}(m^{'})}{n}, 
\end{align*}

which is equivalent to
\begin{align}
\label{inegalité1}
JKL_{\rho,\lambda}^{\otimes_n}(s_0,\hat{s}_{m^{'}}) - \nu_{n}^{\otimes_n} (kl(\bar{s}_m)) &\leq  KL_\lambda^{\otimes_n}(s_0,\bar{s}_m) +\frac{\text{pen}(m) }{n}- \nu_{n}^{\otimes_n} (jkl(\hat{s}_{m^{'}})) \nonumber \\
&+\frac{\eta^{'}+\eta}{n} -\frac{\text{pen}(m^{'})}{n}. 
\end{align}

Mimic the proof as done in Cohen and Le Pennec \cite{Lepennec}, we could obtain that
except on a set of probability less than $ e^{-x_{m^{'}}-x}$, for all $x$, for all $y_{m^{'}} > \sigma_{m^{'}}$, under assumption $(H_m)$, there exists absolute constants $\kappa_0^{'},\kappa_1^{'},\kappa_2^{'}$ such that
\begin{align}
\label{eq nu n}
\frac{- \nu_{n}^{\otimes_n}(jkl(\hat{s}_{m^{'}}))}{y_{m^{'}}^2 + \kappa^{'}_0 d_H^{2 \otimes_n}(s_0,\hat{s}_{m^{'}})} \leq \frac{\kappa_1^{'} \sigma_{m^{'}}}{y_{m^{'}}} + \kappa^{'}_2 \sqrt \frac{x_{m^{'}} +x}{n y_{m^{'}}^2} + \frac{18}{\rho}\frac{x_{m^{'}} + x}{n y_{m^{'}}^2}. 
\end{align}
 
To obtain this inequality we use the hypothesis $(\text{Sep}_m)$ and $(H_m)$. This control is derived from maximal inequalities, described in \cite{MassartStFlour}.

Our purpose is now to control $\nu_n^{\otimes_n}(kl(\bar{s}_m))$.
This is the difference with the theorem of Cohen and Le Pennec: we work with a random subcollection $\mathcal{M}^L$ of $\mathcal{M}$.

By definition of $kl$ and $\nu_n^{\otimes_n}$,
$$\nu_{n}^{\otimes_n} (kl(\bar{s}_m)) = -\frac{1}{n} \sum_{i=1}^n \log \left( \frac{\bar{s}_m(Y_i|X_i)}{s_0 (Y_i|X_i)} \right) + E\left[\frac{1}{n}\sum_{i=1}^n \log \left( \frac{\bar{s}_m(Y_i|X_i)}{s_0 (Y_i|X_i)} \right) \right].$$
We want to apply Bernstein's inequality, which is recalled in appendix.

If we denote by $Z_i$ the random variable $Z_i = -\frac{1}{n} \log \left( \frac{\bar{s}_m (Y_i|X_i)}{s_0(Y_i|X_i)} \right)$, we get $\nu_n^{\otimes_n } (kl (\bar{s}_m)) = \sum_{i=1}^n (Z_i - E(Z_i))$.
We need to control the moments of $Z_i$ to apply Bernstein's inequality.

\begin{lem}
\label{lemma Bernstein}
Let $s_0$ and $\bar{s}_m$ two conditional densities with respect to the Lebesgue measure. Assume that there exists $\tau>0$ such that $\log\left( \left|\left| \frac{s_0}{\bar{s}_m} \right| \right|_\infty \right) \leq \tau$.
Then,
\begin{align*}
\frac{1}{n} \sum_{i=1}^n \int_{\mathbb{R}^q} &\left( \log \left( \frac{s_0(y|x_i)}{\bar{s}_m(y|x_i)} \right) \right)^2 s_0(y|x_i) dy\\
 &\leq \frac{\tau^2}{e^{-\tau}+ \tau -1} KL_\lambda^{\otimes_n} (s_0,\bar{s}_m).
\end{align*}
\end{lem}
We prove this lemma in Appendix \ref{proof lemme Bernstein}.

Because $\frac{\tau^2}{e^{-\tau}+ \tau -1} \underset{\tau \rightarrow \infty}{\sim} \tau$, there exists $A$ such that 
$\frac{\tau^2}{e^{-\tau}+ \tau -1} \leq 2 \tau$ for all $\tau \geq A$.
For $\tau \in ]0,A]$, because this function is continuous and equivalent to $2$ in $0$, there exists $B >0$ such that
$\frac{\tau^2}{e^{-\tau}+ \tau -1} \leq B$.
We obtain that
$\sum_{i=1}^n E(Z_i^2) \leq \frac{1}{n} \delta (1 \vee \tau) KL_\lambda^{\otimes_n}(s_0,\bar{s}_m)$,
where $\delta = 2 \vee B$.

Moreover, for all integers $k \geq 3$,

\begin{align*}
\sum_{i=1}^n E((Z_i)^k_+) & \leq \sum_{i=1}^n \frac{1}{n^k} \int_{\mathbb{R}^q} \left( \log \left( \frac{s_0(y|x_i)}{\bar{s}_m (y|x_i) } \right) \right)_+^k s_0(y|x_i) dy \\
&\leq \frac{n}{n^k} \int_{\mathbb{R}^q} \log \left( \frac{s_0(y|x)}{\bar{s}_m(y|x)} \right) ^{k-2} \log \left( \frac{s_0(y|x)}{\bar{s}_m(y|x)} \right)^2 \mathds{1}_{s_0 \geq \bar{s}_m(y|x)} s_0(y|x) dy \\
&\leq \frac{n}{n^k} \tau^{k-2} \delta (1 \vee \tau) KL^{\otimes_n}_\lambda (s_0,\bar{s}_m).
\end{align*}

Assumptions of Bernstein's inequality are satisfied, with
$$v= \frac{\delta (1\vee \tau) KL^ {\otimes_n}_\lambda(s_0,\bar{s}_m)}{n}, \hspace{1cm} c = \frac{\tau}{n},$$

then, for all $u > 0$, except on a set with probability less than $e^{-u}$,
$$\nu_n^{\otimes_n} (kl (\bar{s}_m)) \leq \sqrt{2vu} + cu.$$
Thus, for all $z>0$, for all $u>0$, except on a set with probability less than $e^{-u}$,
\begin{align}
\label{nu kl sm bar}
\frac{\nu_n^{\otimes_n} (kl (\bar{s}_m))}{z^2 + KL^{\otimes_n}_\lambda (s_0,\bar{s}_m)} \leq \frac{\sqrt{2vu} +cu}{z^2 + KL^{\otimes_n}_\lambda (s_0,\bar{s}_m) } \leq \frac{\sqrt{vu}}{z\sqrt{2 KL^{\otimes_n}_\lambda (s_0,\bar{s}_m)}} + \frac{cu}{z^2}. 
\end{align}

We apply this bound to $u=x+x_m+x_{m'}$.
We get that, except on a set with probability less than $e^{-(x+x_m+x_{m'})}$, using that $a^2+b^2 \geq a^2$, from the inequality \eqref{eq nu n},
$$-\nu_{n}^{\otimes_n } (jkl(\hat{s}_{m'})) \leq (y^2_{m'} + \kappa_{0}' d^{2 \otimes_n}(s_0,\hat{s}_{m'}) )
\left( \frac{\kappa_1' + \kappa_2'}{\theta} + \frac{18}{\theta^2 \rho} \right),$$
and, from the inequality \eqref{nu kl sm bar},
$$\nu_{n}^{\otimes_n } (kl(\bar{s}_m)) \leq (\beta + \beta^2) (z_{m,m'}^2 + KL^{\otimes_n}_\lambda (s,s_m)),$$

where we have chosen
$$y_{m'} = \theta \sqrt{\sigma^2_{m'} + \frac{x_{m'}+x}{n}},$$

with $\theta >1$ to fix later, and
$$z_{m,m'} = \beta^{-1} \sqrt{\left( \frac{v}{2 KL^{\otimes_n}_\lambda(s_0,\bar{s}_m)} + c \right) (x+x_m+x_{m'})},$$
with $\beta >0$ to fix later.

Coming back to the inequality \eqref{inegalité1},
\begin{align*}
JKL_{\rho,\lambda}^{\otimes_n } (s_0,\hat{s}_{m'}) &\leq KL_\lambda^{\otimes_n} (s_0,\bar{s}_m) + \frac{\text{pen}(m)}{n}\\
&+ (y^2_{m'} + \kappa_{0}' d^{2 \otimes_n}(s_0,\hat{s}_{m'}) )
\left( \frac{\kappa_1' + \kappa_2'}{\theta} + \frac{18}{\theta^2 \rho} \right)\\
& + \frac{\eta'+\eta}{n} - \frac{\text{pen}(m')}{n} + (\beta + \beta^2) (z_{m,m'}^2 + KL^{\otimes_n}_\lambda (s_0,\bar{s}_m)).
\end{align*}
Recall that $\bar{s}_m$ is chosen such that
$$KL^{\otimes_n}_\lambda (s_0,\bar{s}_m) \leq \inf_{s_m \in S_m} KL^{\otimes_n}_\lambda (s_0,s_m) + \frac{\delta_{KL}}{n}.$$
Put $\kappa (\beta) = 1 + (\beta + \beta^2)$, and let $\epsilon_1 >0$, we define $\theta_1$ by 
$\kappa_{0}'\left( \frac{\kappa_1' + \kappa_2'}{\theta_1} + \frac{18}{\theta_1^2 \rho} \right) = C_\rho \epsilon_1$ where $C_\rho$ is defined by $C_\rho d_H^{2 \otimes_n } (s_0,\hat{s}_{m'}) \leq JKL_{\rho,\lambda}^{\otimes_n} (s_0,\hat{s}_{m'})$, and put $\kappa_2=\frac{C_\rho \epsilon_1}{\kappa_0}$. We get that
\begin{align*}
(1-\epsilon_1) JKL_{\rho,\lambda}^{\otimes_n} (s_0,\hat{s}_{m'}) &\leq \kappa(\beta) KL^{\otimes_n}_\lambda(s_0,s_m) + \frac{\text{pen}(m)}{n} -\frac{\text{pen}(m')}{n} \\
&+ \kappa(\beta)\frac{\delta_{KL}}{n} + \frac{\eta'+\eta}{n} \\
&+y^2_{m'}  \kappa_2 + (\beta +\beta^2)z_{m,m'}^ 2.
\end{align*}

Since $\tau \leq 1 \vee \tau$, if we choose $\beta$ such that $(\beta + \beta^2) (\delta/2+1)=\alpha \theta_1^{-2} \beta^{-2}$, and putting $\kappa_1 = \alpha \gamma^{-2} (\beta^{-2} +1)$, since $1 \leq 1 \vee \tau$, using the expressions of $y_{m'}$ and $z_{m,m'}$, we get that

\begin{align*}
(1-\epsilon_1) JKL_{\rho,\lambda}^{\otimes_n } (s_0,\hat{s}_{m'}) &\leq \kappa(\beta) KL_\lambda^{\otimes_n}(s_0,s_m) + \frac{\text{pen}(m)}{n} - \frac{\text{pen}(m')}{n} \\
&+\kappa(\beta)\frac{\delta_{KL}}{n} + \frac{\eta' + \eta}{n} \\
& + \kappa_2 \theta_1^2 \left(\sigma_{m'}^2 + \frac{x+x_{m'}}{n}\right) + \kappa_1 (1 \vee \tau) \frac{x+x_m+x_{m'}}{n}\\
&\leq \kappa(\beta)KL^{\otimes n}_\lambda(s_0,s_m)+ \left( \frac{\text{pen}(m)}{n} + \kappa_1(1\vee \tau) \frac{x_m}{n} \right) \\
&+ \left( - \frac{\text{pen}(m')}{n} +\kappa_2 \theta_1^2(\sigma_{m'}^2 + \frac{x_{m'}}{n}) + \kappa_1 (1 \vee \tau) \frac{x_{m'}}{n} \right)\\
&+ \frac{\delta_{KL}}{n} + \frac{\eta' +\eta}{n} + ( \kappa_2 \theta_1^2 + \kappa_1 (1\vee \tau)) \frac{x}{n}.\\
\end{align*}

Now, assume that $\kappa_1 \geq \kappa$ in condition \eqref{penalite}, we get
\begin{align*}
(1-\epsilon_1) JKL^{\otimes_n }_{\rho,\lambda} (s_0,\hat{s}_{m'}) &\leq \kappa(\beta) KL^{\otimes_n}_{\lambda}(s_0,s_m) + 2 \frac{\text{pen}(m)}{n} + \frac{\delta_{KL}}{n} + \frac{\eta+\eta'}{n}\\
&+ (\kappa_2 \theta_1^ 2 + \kappa_1 (1 \vee \tau)) \frac{x}{n}.
\end{align*}
It only remains to sum up the tail bounds over all the possible values of $ m \in \mathcal{M}$ and $m' \in \mathcal{M} (m)$ by taking the union of the different sets of probability less than $e^ {-(x+x_m+x_{m'})}$,
\begin{align*}
\sum_{m \in \mathcal{M} \atop m' \in \mathcal{M}(m)} e^ {-(x+x_m+x_{m'})} &\leq e^{-x} \sum_{(m,m') \in \mathcal{M}\times \mathcal{M}} e^{-(x_m+x_{m'})}\\
&= e^{-x} \left( \sum_{m \in \mathcal{M}} e^{-x_m} \right)^2 = \Sigma^2 e^ {-x} 
\end{align*}
from the assumption $(K)$.

We then have simultaneously for all $m \in \mathcal{M}$, for all $m' \in \mathcal{M}(m)$, except on a set with probability less than $\Sigma^2 e^{-x}$,
\begin{align*}
(1-\epsilon_{1}) JKL_{\rho,\lambda}^{\otimes_n } (s_0,\hat{s}_{m'}) &\leq \kappa(\beta) KL^{\otimes_n}_{\lambda}(s_0,s_m) + 2 \frac{\text{pen}(m)}{n} + \frac{\delta_{KL}}{n} \\
&+ \frac{\eta+\eta'}{n}+ \left(\kappa_2 \theta_1^ 2 + \kappa_1 (1 \vee \tau)\right) \frac{x}{n} .
\end{align*}
It is in particular satisfied for all $m \in \hat{\mathcal{M}}$ and $m' \in \hat{\mathcal{M}}(m)$, and, since $\hat{m} \in \hat{\mathcal{M}}(m)$ for all $m \in \hat{\mathcal{M}}$, we deduce that except on a set with probability less than $\Sigma^2 e^ {-x}$, 
\begin{align*}
JKL_{\rho,\lambda}^{\otimes_n } (s_0,\hat{s}_{m'}) \leq \frac{1}{(1-\epsilon_1) } \times &\left( \inf_{m \in \hat{\mathcal{M}}} \left\{ \kappa(\beta) KL_\lambda^{\otimes_n}(s_0,s_m) + 2 \frac{\text{pen}(m)}{n} \right\} \right. \\
&+ \left. \frac{\delta_{KL}}{n} + \frac{\eta+\eta'}{n}+  \left(\kappa_2 \theta_1^ 2 + \kappa_1 (1 \vee \tau)\right) \frac{x}{n} \right).
\end{align*}
By integrating over all $x>0$, because for any non negative random variable $Z$ and any $a>0$, 
$E(Z) = a \int_{z \geq 0} P(Z > az) dz$, we obtain that
\begin{align*}
&E \left(JKL_{\rho,\lambda}^{\otimes_n } (s_0,\hat{s}_{m'}) - \frac{1}{(1-\epsilon_1) }  \left( \inf_{m \in \hat{\mathcal{M}}} \left\{ \kappa(\beta) KL^{\otimes n}_{\lambda}(s_0,s_m) + 2 \frac{\text{pen}(m)}{n} \right\} \right. \right.\\
&\left. \left.+ \frac{\delta_{KL}}{n} + \frac{\eta+\eta'}{n} \kappa_0 \theta^ 2 \right) \right) \\
&\leq \left(\kappa_2 \theta_1^ 2 + \kappa_1 (1 \vee \tau)\right) \frac{\Sigma^2}{n} .
\end{align*}

As $\delta_{KL}$ can be chosen arbitrary small, this implies that
\begin{align*}
E(JKL^ {\otimes_n} (s_0,\hat{s}_{\hat{m}} )) \leq &\frac{1}{1 - \epsilon_1}  E \left( \inf_{m \in \hat{\mathcal{M}}} \kappa(\beta) KL^ {\otimes_n }_{\lambda} (s_0,s_m)+ \frac{\text{pen}(m)}{n} \right) \\
&+ \frac{\eta + \eta'}{n} + (\kappa_2 \theta_1^ 2 + \kappa_1 (1 \vee \tau) ) \frac{\Sigma^2}{n}\\
\leq &C_1 E \left( \inf_{m \in  \hat{\mathcal{M}}} \inf_{t \in S_m} KL^ {\otimes_n }_{\lambda} (s_0,t)+ \frac{\text{pen}(m)}{n} \right) \\
&+ C_2 ( 1 \vee \tau)\frac{\Sigma^2}{n} + \frac{\eta' + \eta}{n}
\end{align*}
with $C_1= \frac{2}{1-\epsilon_1}$ and $C_2 = \kappa_2 \theta_1^2 + \kappa_1$.

\subsection{Sketch of the proof of the oracle inequality \ref{thm}}
To prove the theorem \ref{thm}, we have to apply the theorem \ref{thm general}. Then, our model has to satisfy all the assumptions.
The assumption $(Sep_m)$ is true when we consider Gaussian densities. 
If $s_0$ is bounded, with compact support, the assumption \eqref{hyp tau} is satisfied. It is also true in others particular cases. We have to look after assumption $(H_m)$ and assumption $(K)$.
Here we present only the main step to prove these assumptions. All the details are in Appendix.
\subsubsection{Assumption $(H_m)$}
We could take $\phi_m (\sigma) = \int_0 ^\sigma \sqrt{ H_{[.]}(\epsilon, S_m,d_H^{\otimes_n})}d\epsilon$ for all $\sigma >0$.
It could be better to consider more local version of the integrated square root entropy, but the global one is enough in this case to define the penalty.
As done in Cohen and Le Pennec \cite{Lepennec}, we could decompose the entropy by
$$H_{[.]} (\epsilon,\mathcal{S}^{B}_{(k,J)}, d_H^{\otimes_n}) \leq H_{[.]} (\epsilon,\Pi_k, d_H^{\otimes_n}) + k H_{[.]} (\epsilon,\mathcal{F}_{J}, d_H^{\otimes_n})$$
where 
\begin{align*}
\mathcal{S}^{\mathcal{B}}_{(k,J)}&= \left\{
\begin{array}{lll}
&y \in \mathbb{R}^q | x \in \mathbb{R}^p \mapsto s_{\theta}(y|x) = \sum_{r=1}^k \pi_r \Phi (y | (\beta_r x) _{|J},\Sigma_r) \\
&\theta = \left\{ \pi_1,\ldots,\pi_k, \beta_1,\ldots,\beta_k,\Sigma_1,\ldots,\Sigma_k \right\} \in \Theta_{(k,J)} \\
&\Theta_{(k,J)}=\Pi_k \times ([-A_\beta,A_\beta]^{|J|})^k \times ([a_\Sigma,A_\Sigma]^q_{+ *})^k
\end{array}
\right\}\\
\Pi_k&=\left\{ (\pi_1,\ldots,\pi_k)\in (0,1)^k ; \sum_{r=1}^k \pi_r =1\right\}\\
\mathcal{F}_{J}&=\left\{ \Phi(.|(\beta X)_{|J},\Sigma) ; \beta \in [a_\beta,A_\beta]^{|J|}, \Sigma = \text{diag}(\Sigma^2_1,\ldots,\Sigma^2_q)\in [a^2_\Sigma,A^2_\Sigma]^q \right\}
 \end{align*}

where $\Phi$ denote the Gaussian density.
\paragraph{Calculus for the proportions}
We could apply a result proved by Wasserman and Genovese in \cite{Wasserman} to bound the entropy for the proportions. We get that
$$H_{[.]} (\epsilon,\Pi_k, d_H^{\otimes_n}) \leq \log \left(k (2 \pi e)^{k/2} \left( \frac{3}{\epsilon}\right)^{k-1} \right).$$
 
\paragraph{Calculus for the Gaussian}
The family
 \begin{equation}
 B_{\epsilon} (\mathcal{F}_{J}) = \left\{
 \begin{array}{lll}
  l(y,x)=(1+\delta)^{-p^2q-3q/4} \Phi(y|\nu_J x,(1+\delta)^{-1/4} B)\\
  u(y,x)=(1+\delta)^{p^2q+3q/4} \Phi(y|\nu_J x,(1+\delta) B) \\
  B=\text{diag}(b_{i(1)}^2,\ldots,b_{i(q)}^2), \text{ with } i \text{ a permutation,}\\
\text{ and } \left\{\begin{array}{lll}
				b_l^2 = (1+\delta)^{1-l/2} A_{\Sigma}^2, l \in \{2,\ldots,R\}\\
				\forall (j,z) \in J^c, \nu_{j,z} = 0 \\
				\forall (j,z) \in J, \nu_{j,z} = \sqrt{ c} \delta A_\Sigma u_{j,z} \\
			\end{array}\right.
 \end{array}
\right\}
\end{equation}
is an $\epsilon$-bracket covering for $\mathcal{F}_{J}$, 
where
$u_{j,z}$ is a net for the mean,
$R$ is the number of parameters needed to recover all the variance set, $\delta=\frac{1}{\sqrt{2}(p^2q+\frac{3}{4}q)} \epsilon$, and $c =\frac{5(1-2^{-1/4})}{8}$.

We obtain that
$$|B_{\epsilon}(\mathcal{F}_{J})| \leq 4 \left( \frac{2 A_\beta}{\sqrt{c} A_\Sigma } \right)^{|J|} \left( \frac{A_\Sigma}{a_\Sigma}+\frac{1}{2} \right) \delta^{1+|J|} ;$$
and then we get
$$H_{[.]} (\epsilon,\mathcal{F}_{J}, d_H^{\otimes_n}) \leq \log \left( 4 \left( \frac{2 A_\beta}{\sqrt{c} A_\Sigma} \right)^{|J|} \left( \frac{A_\Sigma}{a_\Sigma}+\frac{1}{2} \right) \delta^{-1-|J|} \right).$$

\begin{prop}
\label{prop N}
Put $D_{(k,J)}=k(1+|J|)$. For all $\epsilon \in (0,1)$,
$$\mathcal{H}_{[.]} (\epsilon,\mathcal{S}^{B}_{(k,J)}, d_H^{\otimes_n}) \leq \log(C) + D_{(k,J)} \log\left( \frac{1}{\epsilon}\right) ;$$
with
$$C = 4 k (2 \pi e)^{k/2} \left( \frac{2^{5/4} A_\beta}{\sqrt{c} A_\Sigma} \right)^{k |J|} \left( \frac{A_\Sigma}{a_\Sigma}+\frac{1}{2} \right)^k\left(\sqrt{2}q\right)^{k(1+|J|)}.$$
\end{prop}

\paragraph{Determination of a function $\phi$}
We could take
$$\phi_{(k,J)} (\sigma) = \sqrt{D_{(k,J)}} \sigma \left[ B(A_\beta,A_\Sigma,a_\Sigma, q)+\sqrt{\log\left(\frac{1}{\sigma \wedge 1} \right) } \right].$$

This function is non-decreasing, and $\sigma \mapsto \frac{\phi_{(k,J)}(\sigma)}{\sigma}$ is non-increasing.

The root $\sigma_{(k,J)}$ is the solution of $\phi_{(k,J)} (\sigma_{(k,J)}) = \sqrt{n} \sigma_{(k,J)}^2$.
With the expression of $\phi_{(k,J)} $, we get 
$$ \sigma_{(k,J)}^2 = \sqrt{\frac{D_{(k,J)}}{n}} \sigma \left[ B(A_\beta,A_\Sigma,a_\Sigma,q) + \sqrt{\log\left(\frac{1}{\sigma_{(k,J)} \wedge 1} \right) } \right].$$
Nevertheless, we know that $\sigma^*=\sqrt{\frac{D_{(k,J)}}{n}} B(A_\beta,A_\Sigma,a_\Sigma,q)$ minimizes $\sigma_{(k,J)}$: we get
$$ \sigma_{(k,J)}^2 \leq  \frac{D_{(k,J)}}{n} \left[ 2 B^2(A_\beta,A_\Sigma,a_\Sigma,q) + \log\left(\frac{1}{\frac{D_{(k,J)}}{n} B^2(A_\beta,A_\Sigma,a_\Sigma,p,q) \wedge 1} \right)  \right].$$

\subsubsection{Assumption $(K)$}
We want to group models by their dimension.
\begin{lem}
 \label{weightLemma}
 The quantity $\text{card} \{ ( k,J) \in \mathbb{N}^* \times \mathcal{P}([1,p] \times [1,q]), D(k,J)=D\}$ is upper bounded by
 $$ \left\{ \begin{array}{ll}
             &2^{pq} \text{ if } pq \leq D-q^2 \\
             &\left( \frac{epq}{D-q^2}\right)^{D-q^2} \text{ otherwise.}
            \end{array} \right.$$
 \end{lem}
 
  \begin{prop}
  \label{weight}
   Consider the weight family $\{x_{(k,J)} \}_{(k,J)}$ defined by
   $$x_{(k,J)}= D_{(k,J)} \log\left( \frac{4epq}{(D_{(k,J)}-q^2)\wedge pq} \right).$$
   Then we have $\sum_{(k,J)} e^{-x_{(k,J)}} \leq 2$.
  \end{prop}
  \section{Acknowledgment}
I am grateful to Pascal Massart for suggesting me to study this problem, and for stimulating discussions.

 \bibliography{biblio}
\bibliographystyle{plain}

\section{Appendix: technical results}
In this appendix, we give more details for the proofs.
\subsection{Bernstein's lemma}
\begin{lem}[Bernstein's inequality]
 Let $(X_1,\ldots,X_n)$ be independent real valued random variables. Assume that there exists some positive numbers $v$ and $c$ such that 
 $\sum_{i=1}^n E(X_i^2) \leq v$, and, for all integers $k \geq 3$, 
 $\sum_{i=1}^n E((X_i)^k_+ ) \leq \frac{k!}{2} v c^{k-2}$.
 Let $S=\sum_{i=1}^n (X_i - E(X_i))$. Then, for every positive $x$,
 $$P(S \geq \sqrt{2vx} + cx) \leq \exp(-x).$$
\end{lem}

\subsection{Proof of lemma \ref{lemma Bernstein}}
\label{proof lemme Bernstein}
This proof is adapted from the Meynet's thesis, \cite{TheseCaroline}.
First, let give some bounds of functions:
\begin{lem}
 Let $\tau >0$. For all $x>0$, consider
 $$f(x)= x \log(x)^2, \hspace{1cm} h(x)=x \log(x) -x+1, \hspace{1cm} \phi(x)= e^x -x-1.$$
 Then, for all $0<x<e^\tau$, we get
 $$f(x) \leq \frac{\tau^2}{\phi(-\tau)} h(x).$$
\end{lem}
To prove this, we have to show that $y \mapsto \frac{\phi(y)}{y^2}$ is non-decreasing. We omit the proof here.

We want to apply this inequality, in order to derive the lemma \ref{lemma Bernstein}.
As $\log\left(\left|\left|\frac{s}{\bar{s}_m}\right|\right|_\infty\right) \leq \tau$,
$$\left|\left|\frac{s_0}{\bar{s}_m} \right|\right| _\infty \leq e^{\tau};$$
and we could apply the previous inequality to $\frac{s_0}{\bar{s}_m}$.
Indeed,
$$f\left(\frac{s_0}{\bar{s}_m}\right) \leq \frac{\tau^2}{\phi(-\tau)} h \left(\frac{s_0}{\bar{s}_m} \right).$$
Integrating with respect to the density $\bar{s}_m$, we get that
\begin{align*}
&\int \frac{s_0(y|.)}{\bar{s}_m(y|.)} \log\left( \frac{s_0(y|.)}{\bar{s}_m(y|.)} \right)^2 \bar{s}_m(y|.) dy \\
\leq &\int \frac{\tau^2}{e^{-\tau}-\tau-1} \left( \frac{s_0(y|.)}{\bar{s}_m(y|.)} \log \frac{s_0(y|.)}{\bar{s}_m(y|.)} - \frac{s_0(y|.)}{\bar{s}_m(y|.)}+1 \right) \bar{s}_m(y|.) dy\\
\Longleftrightarrow &\frac{1}{n} \sum_{i=1}^n \int s_0(y|x_i) \log\left( \frac{s_0(y|x_i)}{\bar{s}_m(y|x_i)} \right)^2 dy \\
\leq & \frac{\tau^2}{e^{-\tau}-\tau-1} \frac{1}{n} \sum_{i=1}^n \int  s_0(y|x_i) \log \frac{s_0(y|x_i)}{\bar{s}_m(y|x_i)} dy.
 \end{align*}
This conclude the proof.

\subsection{Determination of a net for the mean and the variance}
\begin{itemize}
 \item \textbf{Step 1: construction of a net for the variance}
 
 Let $\epsilon \in ]0,1]$, and $\delta = \frac{1}{\sqrt{2}(p^2q+\frac{3}{4}q} \epsilon$.
Let $b_j^2 = (1+\delta)^{1-\frac{j}{2}} A_\Sigma^2$.
For $2 \leq j \leq R$, we have \\ \mbox{$[a_\Sigma, A_\Sigma]= [b_R,b_{R-1}] \bigcup \ldots \bigcup [b_3,b_2]$}, where $R$ is chosen to recover everything.
We want that 
\begin{align*}
\phantom{\Leftrightarrow }\hspace{1cm}& a_\Sigma^2 = (1+\delta)^{1-R/2} A_\Sigma^2\\
\Leftrightarrow \hspace{1cm}& 2 \log \frac{a_\Sigma}{A_\Sigma} = \left( 1-\frac{R}{2} \right) \log(1+\delta)\\
\Leftrightarrow \hspace{1cm}& R = \frac{4 \log( \frac{A_\Sigma}{a_\Sigma}\sqrt{1+\delta} )}{\log (1+\delta)}.
\end{align*}

We want $R$ to be an integer, then $R = \left\lceil \frac{4 \log( \frac{A_\Sigma}{a_\Sigma}\sqrt{1+\delta} )}{\log (1+\delta)}  \right\rceil$.
We get a net for the variance.
We could let $B=\text{diag}(b_{i(1)}^2,\ldots, b_{i(q)}^2)$, close to $\Sigma$ (and deterministic, independent of the values of $\Sigma$), where $i$ is a permutation such that $b_{i(z)+1} \leq \Sigma_z \leq b_{i(z)}$ for all $z \in [1,q]$.
Remember that $\frac{b_{j+1}^2}{b_j^2}= \frac{1}{\sqrt{1+\delta}}$, and that if $\Sigma$ is fixed, $\Sigma= \text{diag} (\Sigma_1^2,\ldots,\Sigma_q^2)$.
 
 \item \textbf{Step 2: construction of a net for the mean vectors}

We select only the active variables detected by the Lasso.
$$J = \left\{ (j,z) \in [1,p] \times [1,q] | \hat{\beta}^{\text{Lasso}}_{j,z} \neq 0\right\}.$$

 Let $f = \Phi(.|\beta x, \Sigma) \in \mathcal{F}_J$.

 \begin{itemize}
 
 \item \textbf{Definition of the brackets}
 
 Define the bracket by the functions $l$ and $u$:
 \begin{align*}
  l(y,x)&=(1+\delta)^{-p^2q-3q/4} \Phi(y|\nu_J x,(1+\delta)^{-1/4} B) ;\\
  u(y,x)&=(1+\delta)^{p^2q+3q/4} \Phi(y|\nu_J x,(1+\delta) B).
 \end{align*}
We have chosen $i$ such that $b_{i(z)+1}^2 \leq \Sigma_z^2 \leq b_{i(z)}^2$ for all $1\leq z \leq q$.

We need to define $\nu$ such that $[l,u]$ is an $\epsilon$-bracket for $f$.

\item \textbf{ Proof that $[l,u]$ is an $\epsilon$-bracket for $f$}

We are looking for a condition on $\nu_J$ to have $ \frac{f}{u} \leq 1$ and $\frac{l}{f} \leq 1$.

We will use the following lemma to compute these ratios.

\begin{lem}
 Let $\Phi(.|\mu_1,\Sigma_1)$ and $\Phi(.|\mu_2,\Sigma_2)$ be two Gaussian densities. If their variance matrices are assumed to be diagonal, with $\Sigma_a=\text{diag}(S_{a1}^2,\ldots,S_{aq}^2)$ for $a \in \{1,2\}$, such that $S^2_{2z}>S_{1z}^2>0$ for all $z \in \{1,\ldots,q\}$,
 then, for all $x \in \mathbb{R}^q$,
\begin{align*}
\frac{\Phi(x|\mu_1,\Sigma_1)}{\Phi(x|\mu_1,\Sigma_1)} \leq \prod_{z=1}^q \frac{\sqrt{\Sigma_{2z}}}{\sqrt{\Sigma_{1z}}} e^{\frac{1}{2} (\mu_1-\mu_2)^t \text{diag} \left( \frac{1}{\Sigma_{21}-\Sigma_{11}},\ldots,\frac{1}{\Sigma_{2q}-\Sigma_{1q}}\right) (\mu_1-\mu_2) }.
\end{align*}

 \end{lem}
 
 For the ratio $\frac{f}{u}$ we get:
\begin{align}
\label{fu}
\frac{f(y|x)}{u(y,x)} = &\frac{1}{(1+\delta)^{p^2q+3q/4}} \frac{\Phi(y|\beta x,\Sigma)}{\Phi(y|\nu_J x,(1+\delta)B)} \\
\leq &\frac{1}{(1+\delta)^{p^2q+3q/4}} \prod_{z=1}^q \frac{b_z}{\Sigma_z} (1+\delta)^{q/2} \nonumber\\
&\times e^{\frac{1}{2} (\beta x - \nu_J x)^t ((1+\delta)B-\Sigma)^{-1} (\beta x - \nu_J x)} \nonumber \\
\leq &(1+\delta)^{p^2q-q/4} (1+\delta)^{q/4} e^{ \frac{1}{2} (\beta x - \nu_J x)^t (\delta B)^{-1} (\beta x - \nu_J x)} \nonumber \\
\leq &(1+\delta)^{p^2q} e^{ \frac{1}{2 \delta} (\beta x - \nu_J x)^t B^{-1} (\beta x - \nu_J x)}.\nonumber
\end{align}

For the ratio $\frac{l}{f}$ we get:
\begin{align}
\label{lf}
\frac{l(y,x)}{f(y|x)} = &\frac{1}{(1+\delta)^{p^2q+3q/4}} \frac{\Phi(y|\nu_J  x,(1+\delta)^{-1/4} B)}{\Phi(y|\beta x,\Sigma)} \\
\leq &\frac{1}{(1+\delta)^{p^2q+3q/4}} \prod_{z=1}^q \frac{\Sigma_z}{b_z} (1+\delta)^{q/8}\nonumber\\
&\times e^{ \frac{1}{2} (\beta x - \nu_J x)^t (\Sigma-B)^{-1} (\beta x - \nu_J x) } \nonumber \\
\leq &(1+\delta)^{-p^2q-3q/8} (1+\delta)^{q/4} \nonumber\\
&\times e^{ \frac{1}{2} (\beta x - \nu_J x)^t ((1-(1+\delta)^{-1/4}) B)^{-1} (\beta x - \nu_J x)} \nonumber \\
\leq &(1+\delta)^{-p^2q-3q/8} e^{ \frac{1}{2(1-(1+\delta)^{-1/4})} (\beta x - \nu_J x)^t B^{-1} (\beta x - \nu_J x)}.\nonumber
\end{align}

We want to bound the ratios \eqref{fu} and \eqref{lf} by $1$.
Put $c=\frac{5 (1-2^{-1/4})}{8}$, and develop these calculus.
A necessary condition to obtain this bound is 
$$||\beta x - \nu_Jx||_2^2 \leq pq \delta^2 (1-2^{-1/4})A_\Sigma^2 .$$
Indeed, we want
\begin{align*}
 (1+\delta)^{-p^2q-3q/8} e^{\frac{1}{2(1-(1+\delta)^{-1/4})} (\beta x - \nu_J x)^t B^{-1} (\beta x - \nu_J x)} &\leq 1 \\
 (1+\delta)^{-p^2q} e^{ \frac{1}{2 \delta A_{\Sigma}} (\beta x - \nu_J x)^t B^{-1} (\beta x - \nu_J x) }& \leq 1 ;
\end{align*}
which is equivalent to
\begin{align*}
 ||\beta x-\nu_J x||^2_2 &\leq p^2q \frac{\delta^2}{2} A_\Sigma^2 ;\\
 ||\beta x-\nu_J x||^2_2 &\leq (p^2q+\frac{3}{4}q) \delta^2 (1-2^{-1/4})A_\Sigma.
\end{align*}
As $||\beta x-\nu_J x||_2^2 \leq p||\beta-\nu_J||_2^2 ||x||_\infty$, and $X \in [0,1]^p$, we need to get \\
$||\beta-\nu_J||_2^2 \leq p q \delta^2 (1-2^{-1/4})A_\Sigma^2$ to have the wanted bound.
Put
$$ U:= \mathbb{Z} \cap \left[ \left\lfloor \frac{-A_\beta}{\sqrt{ c } \delta  A_\Sigma} \right\rfloor, \left\lfloor \frac{A_\beta}{\sqrt{ c}\delta A_\Sigma} \right\rfloor \right].$$
For all $j \in J$, choose
$$u_{j,z}= \underset{v_{j,z} \in U}{\operatorname{argmin}} \left|\beta_{j,z} - \sqrt{c} \delta A_\Sigma v_{j,z}\right|.$$
Define $\nu$ by

\begin{align*}
\text{for all } (j,z) &\in J^c, \nu_{j,z} = 0 ; \\
\text{for all } (j,z) &\in J\phantom{^c}, \nu_{j,z} = \sqrt{c} \delta  A_\Sigma u_{j,z}. 
\end{align*}

Then, we get a net for the mean vectors.
 
\item \textbf{Proof that $[l,u]$ is an $\epsilon$-bracket}

We will work with the Hellinger distance.

\begin{align*}
 d_H^2(l,u)&= \frac{1}{2} \int_{\mathbb{R}^q} (\sqrt{l}-\sqrt{u})^2 d\lambda\\
 &=\frac{1}{2} \int_{\mathbb{R}^q} l+u - 2\sqrt{lu} d\lambda\\
 &=\frac{1}{2} \left[(1+\delta)^{-p^2q-3q/4} + (1+\delta)^{p^2q+3q/4} \right] - \int_{\mathbb{R}^q} \sqrt{\Phi_l \Phi_u} d\lambda\\
 &=\frac{1}{2} \left[(1+\delta)^{-p^2q-3q/4} + (1+\delta)^{p^2q+3q/4} \right] \\
 &- \left( \prod_{z=1}^q \frac{2 b_{i(z)+1} b_{i(z)} (1+\delta)^{1/2} (1+\delta)^{-1/8}}{(1+\delta)b_{i(z)+1}^2+ (1+\delta)^{-1/4} b_{i(z)}}^2 \right)^{1/2} * 1 . \\
\end{align*}
We have used the following lemma:
\begin{lem}
The Hellinger distance of two Gaussian densities with diagonal variance matrices is given by the following expression:
\begin{align*}
 &d_H^2(\Phi(.|\mu_1,\Sigma_1),\Phi(.|\mu_2,\Sigma_2))\\
 = & 2-2 \left(\prod_{q_1=1}^{q}\frac{2 \Sigma_{1 q_1}\Sigma_{2 q_1}}{\Sigma_{1q_1}^2 + \Sigma_{2q_1}^2}\right)^{1/2} \\
 &\times \exp\left\{-\frac{1}{4} (\mu_1-\mu_2)^t \text{diag} \left( \left(\frac{1}{\Sigma_{1 q_1}^2+\Sigma_{2q_1}^2} \right)_{q_1=1,\ldots,q} \right) (\mu_1-\mu_2)\right\}
\end{align*}
 \end{lem}

As $b_{i(z)+1}^2 = (1+\delta)^{-1/2} b_{i(z)}^2$, we get that
\begin{align*}
2 \frac{ (1+\delta)^{3/8} b_{i(z)}^2}{b_{i(z)+1}^2 \left[ (1+\delta)^{-1/4} + (1+\delta)^{1/2}(1+\delta) \right] } &= 2 \frac{(1+\delta)^{5/8}}{(1+\delta)^{-1/4}  + (1+\delta)^{3/2}} \\
&= \frac{2}{(1+\delta)^{-7/8}  + (1+\delta)^{7/8}} .
\end{align*}
Then

\begin{align*}
d_H^2(l,u)=&\frac{1}{2} \left[(1+\delta)^{-(p^2q+3q/4)} + (1+\delta)^{p^2q+3q/4} \right] \\
&- \left(\frac{2}{(1+\delta)^{-7/8} + (1+\delta)^{7/8}} \right)^{q/2}\\
d_H^2(l,u)=&\cosh((p^2q+3q/4) \log(1+\delta)) - 2 \cosh(7/8 \log(1+\delta))^{-q/2}\\
=&\cosh((p^2q+3q/4) \log(1+\delta)) -1 + 1\\
&- 2^{-q/2} \cosh(7/8 \log(1+\delta))^{-q/2}.\\
\end{align*}

We want to apply the Taylor formula to $f(x)=\cosh(x)-1$ to obtain an upper bound, and to $g(x)=1-2^{-q/2} \cosh(x) ^{-q/2}$.
Indeed, there exists $c$ such that, on the good interval, $f(x) \leq \cosh(c) \frac{x^2}{2}$ and $g(x) \leq q^2 \frac{x^2}{2}$.
Then, and because $\log(1+\delta) \leq \delta$, 

\begin{align*}
d_H^2(l,u)& \leq \cosh((p^2q+3q/4) \log(1+\delta)) - 2 \cosh(7/8 \log(1+\delta))^{-q/2}\\
&\leq (p^2q+3q/4)^2 \delta^2 \left( \cosh(\alpha) + \frac{49}{128} \right)\\
&\leq 2 (p^2q+3q/4)^2 \delta^2 \leq \epsilon^2.
\end{align*}
where $\epsilon \geq \sqrt{2} (p^2 q +\frac{3}{4}q)\delta$.
\end{itemize}
  
 \item \textbf{Step 3: Upper bound of the number of $\epsilon$-brackets for $\mathcal{F}_{J}$.}
 
 From step $1$ and step $2$, the family
 \begin{equation}
 B_{\epsilon} (\mathcal{F}_{J}) = \left\{
 \begin{array}{lll}
  l(y,x)=(1+\delta)^{-(p^2q+3q/4)} \Phi(y|\nu_J x,(1+\delta)^{-1/4} B)\\
  u(y,x)=(1+\delta)^{p^2q+3q/4} \Phi(y|\nu_J x,(1+\delta) B) \\
  B=\text{diag}(b_{i(1)},\ldots, b_{i(q)})\text{ where i is a permutation}\\
\text{ with } \left\{\begin{array}{lll}
				b_{i(z)}^2 = (1+\delta)^{1-i(z)/2} A_{\Sigma}^2 \text{ for all } z\in \{1,\ldots,q\}\\
				\forall (j,z) \in J^c, \nu_{j,z} = 0 \\
				\forall (j,z) \in J, \nu_{j,z} = \sqrt{ c} \delta A_\Sigma u_{j,z} \\
			\end{array}\right.
 \end{array}
\right\}
\end{equation}
is an $\epsilon$-bracket for $\mathcal{F}_{J}$.
Therefore, an upper bound of the number of $\epsilon$-brackets necessary to cover $\mathcal{F}_{J}$ is deduced from an upper bound of the cardinal of $B_{\epsilon}(\mathcal{F}_{J})$.
 \begin{align*}
 |B_{\epsilon}(\mathcal{F}_{J})|&\leq \sum_{l=2}^R \prod_{(j,z) \in J} \left( \frac{2 A_\beta}{\sqrt{c} \delta A_\Sigma} \right) \\
 &\leq \left( \frac{2 A_\beta}{\sqrt{c}  \delta A_\Sigma} \right)^{|J|} \sum_{l=2}^R 1\\
 &\leq \left( \frac{2 A_\beta}{\sqrt{c}  \delta A_\Sigma} \right)^{|J|} (R-1).
 \end{align*}
\end{itemize}
But $R\leq \frac{4 \left(\frac{A_\Sigma}{a_\Sigma} +1/2 \right)}{\delta}$, then we get
$$|B_{\epsilon}(\mathcal{F}_{J})| \leq 4 \left( \frac{2 A_\beta}{\sqrt{c} A_\Sigma} \right)^{|J|} \left( \frac{A_\Sigma}{a_\Sigma}+\frac{1}{2} \right)\delta^{-1-|J|}$$

\subsection{Calculus for the function $\phi$}
\label{annexe calcul phi}
From the proposition \ref{prop N}, we obtain, for all $\xi>0$,
$$\int_{0}^\xi \sqrt{\mathcal{H}_{[.]} (\epsilon,\mathcal{S}^{B}_{(k,J)}, d_H^{\otimes_n})} d\epsilon \leq \xi \sqrt{\log(C)} + \sqrt{D_{(k,J)}} \int_{0}^{\xi\wedge1} \sqrt{\log \left(\frac{1}{\epsilon}\right)} d\epsilon$$
we need to control $\int_0^\xi \sqrt{ \log\left( \frac{1}{\epsilon} \right) } d\epsilon$, which is done in Maugis and Meynet (\cite{Maugis}).
\begin{lem}
For all $\xi >0$,
$$\int_{0}^{\xi} \sqrt{\log \left(\frac{1}{\epsilon}\right) }d\epsilon \leq \xi \left[ \sqrt{\pi} + \sqrt{\log\left( \frac{1}{\xi} \right) } \right].$$
\end{lem}
Then
\begin{align*}
\int_{0}^\xi \sqrt{\mathcal{H}_{[.]} (\epsilon,\mathcal{S}^{B}_{(k,J)}, d_H^{\otimes_n})} d\epsilon \leq &\xi \sqrt{\log(C)} \\
&+ \sqrt{D_{(k,J)}} (\xi \wedge 1) \left[ \sqrt{\pi} + \sqrt{\log\left( \frac{1}{\xi\wedge 1} \right) } \right] \\
\leq &\xi \sqrt{D_{(k,J)}} \left[ \sqrt{\frac{\log(C)}{D_{(k,J)}}} + \sqrt{\pi} + \sqrt{\log \left( \frac{1}{\xi\wedge 1}\right) } \right]
\end{align*}

But
\begin{align*}
 \log(C) &\leq \log(4) + \log(k) + \frac{k}{2} \log(2 \pi e)  \\
 &+ k |J| \log \left( \frac{2^{5/4} A_\beta}{\sqrt{c}A_\Sigma} \right) +  k \log \left(\frac{A_\Sigma}{a_\Sigma} + \frac{1}{2} \right) + D_{(k,J)}\log(\sqrt{2}q) + \log(k) \\
 &\leq D_{(k,J)} \left[ \log(4) + \log(\sqrt{2 \pi e}) +\log(\sqrt{2}q) \right. \\
 & \left.+\log \left(\frac{A_\Sigma}{a_\Sigma} + \frac{1}{2} \right) + \log(k) + \log \left( \frac{2^{5/4} A_\beta}{\sqrt{c}A_\Sigma} \right) \right]\\
&\leq  D_{(k,J)}\left[\log(q) + \log \left(\frac{A_\beta}{A_\Sigma}\left(\frac{A_\Sigma}{a_\Sigma} + \frac{1}{2} \right) \right) + \log\left( \sqrt{\pi e} \frac{2^{5/4} 8}{\sqrt{c}} e \right) \right].
 \end{align*}

 Then
 \begin{align*}
&\int_{0}^\xi \sqrt{\mathcal{H}_{[.]} (\epsilon,\mathcal{S}^{B}_{(k,J)}, d_H^{\otimes_n})} d\epsilon \\
\leq &\xi \sqrt{D_{(k,J)}} \left[ \sqrt{\log(q) + \log \left(\frac{A_\beta}{A_\Sigma}\left(\frac{A_\Sigma}{A_\Sigma} + \frac{1}{2} \right) \right) + \log\left( \sqrt{\pi e} \frac{2^{5/4} 8 }{\sqrt{c}} e \right)} \right.\\
&\hspace{1.5cm}+  \left. \sqrt{\pi} + \sqrt{\log \left(\frac{1}{\xi\wedge 1} \right)} \right]\\
\leq& \xi \sqrt{D_{(k,J)}} \left[ \sqrt{\log(q)} + \sqrt{ \log \left(\frac{A_\beta}{A_\Sigma}\left(\frac{A_\Sigma}{a_\Sigma} + \frac{1}{2} \right) \right)}\right.\\
&\left.\hspace{1.5cm}+ a +  \sqrt{\log \left(\frac{1}{\xi\wedge 1} \right)}  \right]\\
\leq& \xi \sqrt{D_{(k,J)}}  \left[ B(A_\beta,A_\Sigma,a_\Sigma,q)+\sqrt{\log\left(\frac{1}{\xi \wedge 1} \right) } \right] ;
\end{align*}
with
$$B(A_\beta,A_\Sigma,a_\Sigma,q)=\sqrt{\log(q)} + \sqrt{ \log \left(\frac{A_\beta}{A_\Sigma}\left(\frac{A_\Sigma}{a_\Sigma} + \frac{1}{2} \right) \right)} + a ;$$
and $a=\sqrt{\pi} + \sqrt{\log(\sqrt{\pi e} 2^{5/4}\frac{8e}{\sqrt{c}})}$.

\subsection{Proof of the proposition \ref{weight}}
We are interested by $\sum_{(k,J) \in \mathcal{M}} e^{-x_{(k,J)}}$. 
Considering 
$$x_{(k,J)}=D_{(k,J)} \log \left( \frac{4epq}{(D_{(k,J)}-q^2) \wedge pq} \right),$$
we could group models by their dimension to compute this sum. Denote by $C_D$ the cardinal of models of dimension $D$.
 \begin{align*}
  &\sum_{(k,J) \in \mathbb{N}^* \times [1,p]\times [1,q]} e^{- D_{(k,J)} \log\left( \frac{4epq}{(D_{(k,J)}-q^2)\wedge pq} \right)} = \sum_{D \geq 1} C_D e^{- D \log\left( \frac{4epq}{(D-q^2)\wedge pq} \right)}\\
  &= \sum_{D=1}^{pq+q^2} e^{- D \log\left( \frac{4epq}{(D-q^2)} \right)}\left( \frac{epq}{D-q^2}\right)^{D-q^2} + \sum_{D=pq+q^2+1}^{+\infty} e^{- D \log\left( \frac{4epq}{pq} \right)}2^{pq} \\
  &= \sum_{D=1}^{pq+q^2}  4^{-D} \left( \frac{epq}{D-q^2}\right)^{-q^2}  + \sum_{D= pq+q^2+1}^{+\infty} e^{-D(\log(4)+1) + pq \log(2)} \\
  &\leq \sum_{D=1}^{pq+q^2}  2^{-D}  + \sum_{D= pq+q^2+1}^{+\infty} 2^{-D} = 2.
 \end{align*}
\subsection{Proof of the lemma \ref{weightLemma}}

  We know that $D_{(k,J)}=k-1+|J|k+kq^2$.
  Then, 
 \begin{align*}
C_D&=\text{card} \{ ( k,J) \in \mathbb{N}^* \times \mathcal{P}([1,p] \times [1,q]), D(k,J)=D\} \\
&\leq \sum_{k \in \mathbb{N}^*} \sum_{(j,z) \in [1,p]\times[1,q]} \binom{pq}{|J|} \mathds{1}_{k(|J|+q^2+1)-1=D}\\
  &\leq \sum_{|J| \in \mathbb{N}^*} \binom{pq}{|J|} \mathds{1}_{|J| \leq pq \wedge (D-q^2)}.
 \end{align*}
 If $pq < D-q^2$,
 $$\sum_{|J|>0} \binom{pq}{|J|} \mathds{1}_{|J| \leq pq \wedge (D-q^2)} = 2^{pq}.$$
 Otherwise, according to the proposition $2.5$ in Massart (\cite{MassartStFlour}),
 $$\sum_{|J|>0} \binom{pq}{|J|} \mathds{1}_{|J| \leq pq \wedge (D-q^2)} \leq f(D-q^2)$$
 where $f(x)=\left( \frac{epq}{x} \right)^x$ is an increasing function on $[1,pq]$.
 As $pq$ is an integer, we get the result.

\end{document}